\documentclass[11pt,a4paper]{article}
\usepackage{amsmath, amsfonts, amsthm, amssymb, graphicx, color}
\usepackage{enumerate, hyperref}

\begin{document}

\title
{4-regular 4-connected Hamiltonian graphs with few Hamiltonian cycles}
\author{
{\sc Carsten THOMASSEN\footnote{Department of Applied Mathematics and Computer Science, Technical University of Denmark, DK-2800 Lyngby, Denmark. E-mail address: ctho@dtu.dk}} \ and {\sc Carol T. ZAMFIRESCU\footnote{Department of Applied Mathematics, Computer Science and Statistics, Ghent University, Krijgslaan 281 - S9, 9000 Ghent, Belgium  and Department of Mathematics, Babe\c{s}-Bolyai University, Cluj-Napoca, Roumania. E-mail address: czamfirescu@gmail.com}}}
\date{}

\maketitle
\begin{center}
\vspace{2mm}
\begin{minipage}{125mm}
{\bf Abstract.}
We prove that there exists an infinite family of 4-regular 4-connected Hamiltonian graphs with a bounded number of Hamiltonian cycles. We do not know if there exists such a family of 5-regular 5-connected Hamiltonian graphs.
\smallskip


\smallskip


\end{minipage}
\end{center}

\vspace{15mm}

\section{Introduction}

There is a variety of 3-regular 3-connected graphs with no Hamiltonian cycles. Much less is known about 4-regular 4-connected graphs. Thus the Petersen graph (on 10 vertices) is the smallest 3-regular 3-connected non-Hamiltonian graph whereas it was an open problem of Nash-Williams if there exists a 4-regular 4-connected non-Hamiltonian graph until Meredith~\cite{Me73} gave an infinite family, the smallest of which has 70 vertices, see~\cite[p.~239]{BM76}. Tait's conjecture that every 3-regular 3-connected graph has a Hamiltonian cycle was open from 1880 till Tutte~\cite{Tu46} in 1946 found a counterexample, see~\cite[p.~161]{BM76}. By Steinitz' theorem, the Tutte graph (and subsequently many others) also show that there are infinitely many 3-regular 3-polyhedral graphs which are non-Hamiltonian, whereas it is a longstanding conjecture of Barnette that every 4-regular 4-polyhedral graph has a Hamiltonian cycle~\cite[p.~1145]{G70} (see also \cite[p.~389a]{G03}
and~\cite[p.~375]{G67}). There are infinitely many 3-regular hypohamiltonian graphs whereas it is a longstanding open problem if there exists a hypohamiltonian graph of minimum degree at least 4, see~\cite{Th78-2}. Starting with the complete graph on 4~vertices and successively replacing vertices by triangles we obtain an infinite family of graphs with precisely three Hamiltonian cycles. Cantoni's conjecture says that these are precisely the planar 3-regular graphs with exactly three Hamiltonian cycles, see~\cite{Tu76}.

Recently, Haythorpe~\cite{Ha18} conjectured that 4-regular graphs behave differently from the 3-regular graphs, also in this respect, in that the number of Hamiltonian cycles increases as a function of the number of vertices. The purpose of this note is to answer this in the negative. So in this respect, the 4-regular 4-connected graphs behave in a similar way as the 3-regular 3-connected graphs. We do not know if this also holds for the 5-regular 5-connected graphs.

Andrew Thomason~\cite{Th78} proved that every Hamiltonian graph whose vertices are all of odd degree has a second Hamiltonian cycle. Thomason's theorem was inspired by (and extends) Smith's result stating that every cubic graph has an even number of Hamiltonian cycles through each edge, see~\cite{Tu46}. Sheehan~\cite{Sh75} conjectured that the same holds for 4-regular graphs: if they are Hamiltonian, then they contain at least two Hamiltonian cycles. In~\cite{GMZ20} this was shown to hold up to order~21. If true in general, this would imply that, for every natural number $k \ge 3$, every $k$-regular Hamiltonian graph has a second Hamiltonian cycle. The latter statement was verified in~\cite{Th98} for all $k>72$ and subsequently in~\cite{HSV07} for all $k>20$.

For restricted classes of 4-regular 4-connected graphs, it is still possible that the number of Hamiltonian cycles in a Hamiltonian graph must increase (perhaps even exponentially) as a function of the number of vertices. This may be true for planar graphs where it is known that the number of Hamiltonian cycles increases at least as a linear function,~\cite{BV21}. And it may hold for bipartite graphs where it is known that the number of Hamiltonian cycles in a Hamiltonian $k$-regular graph increases more than exponentially as a function of $k$, \cite{Th96}.

\section{4-regular 4-connected graphs with a bounded number of Hamiltonian cycles}

As mentioned earlier, the Meredith graph is a 4-regular 4-connected non-Hamiltonian graph.

\vspace{3mm}

Figure 1 indicates an infinite family of 4-regular Hamiltonian graphs with a bounded number of Hamiltonian cycles.

\begin{center}
\includegraphics[height=45mm]{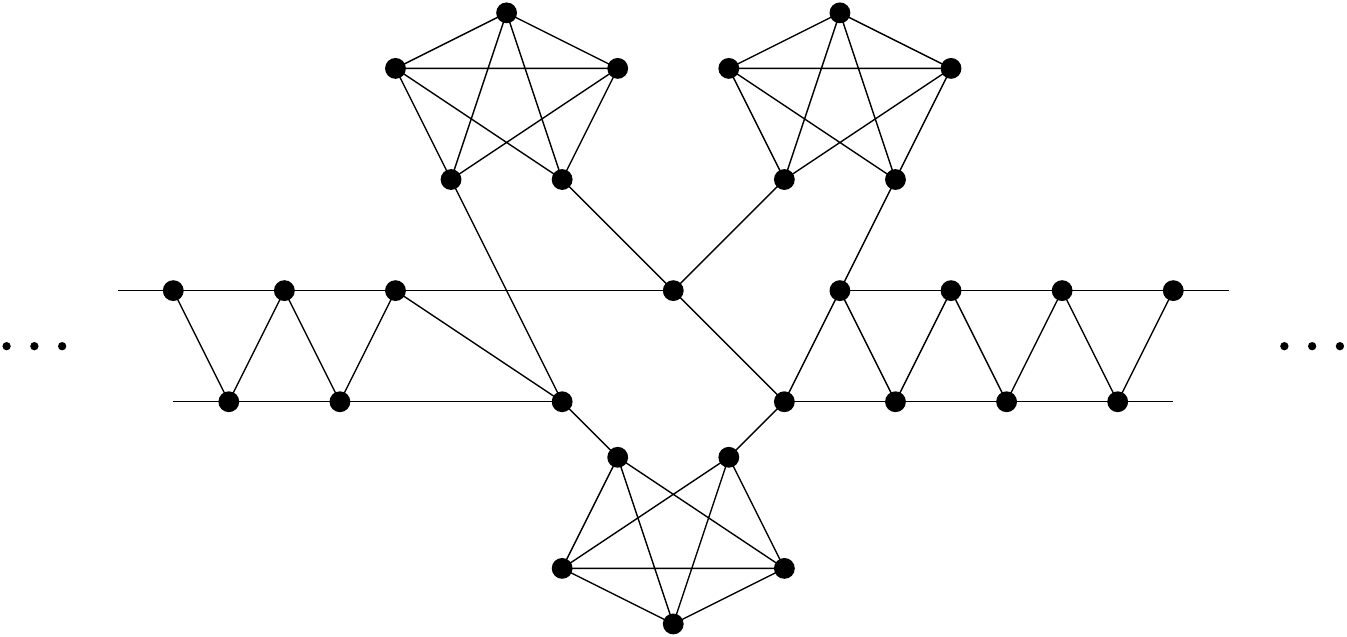}\\
Fig.~1: 4-regular graphs, each having exactly 216~Hamiltonian cycles.\\
(The left-most and right-most part of the graph are to be connected in the obvious way.)
\end{center}

We shall now describe 4-regular 4-connected graphs with a bounded (positive) number of Hamiltonian cycles.

\bigskip

\noindent \textbf{Theorem.} \emph{There exists a constant $c > 0$ such that there are infinitely many $4$-regular $4$-connected graphs, each containing exactly $c$ Hamiltonian cycles.}

\bigskip

\noindent \emph{Proof.} Assume $G$ is a 4-regular 4-edge-connected graph containing a path $abcd$ such that\\[1mm]
(i) $G$ has no Hamiltonian cycle;\\
(ii) $G$ has a $2$-factor consisting of two cycles $C$ and $C'$ such that $C$ contains $ab$ and $C'$ contains $cd$; and\\
(iii) $G - bc$ has no Hamiltonian path joining two of $a,b,c,d$. If $v$ is any vertex in $\{a,b,c,d\}$, then $G - v - bc$ has no Hamiltonian path joining two of $a,b,c,d$.
\smallskip

We construct the graph $H_G$ as follows. Let $G'$ be a copy of $G$ and $\ell \ge 1$ a natural number. Take the disjoint union of $G,G'$; denote for a vertex $v$ in $G$ its copy in $G'$ by $v'$; delete the edges $ab,bc,cd,a'b',b'c',c'd'$; and add four pairwise disjoint paths $P_a, P_b, P_c, P_d$ joining $a,c'$ and $b,b'$ and $c,a'$ and $d,d'$, respectively, where $P_a$ and $P_d$ have length $\ell + 1$, while $P_b$ and $P_c$ have length $\ell$. Thereafter, add a zig-zag path between $P_a$ and $P_b$, as well as a zig-zag path between $P_c$ and $P_d$. The construction of $H_G$ is illustrated in Fig.~2. Every Hamiltonian cycle in $H_G$ will contain the paths $P_a, P_b, P_c, P_d$. Note that the number of Hamiltonian cycles of $H_G$ is independent of $\ell$.

\begin{center}
\includegraphics[height=65mm]{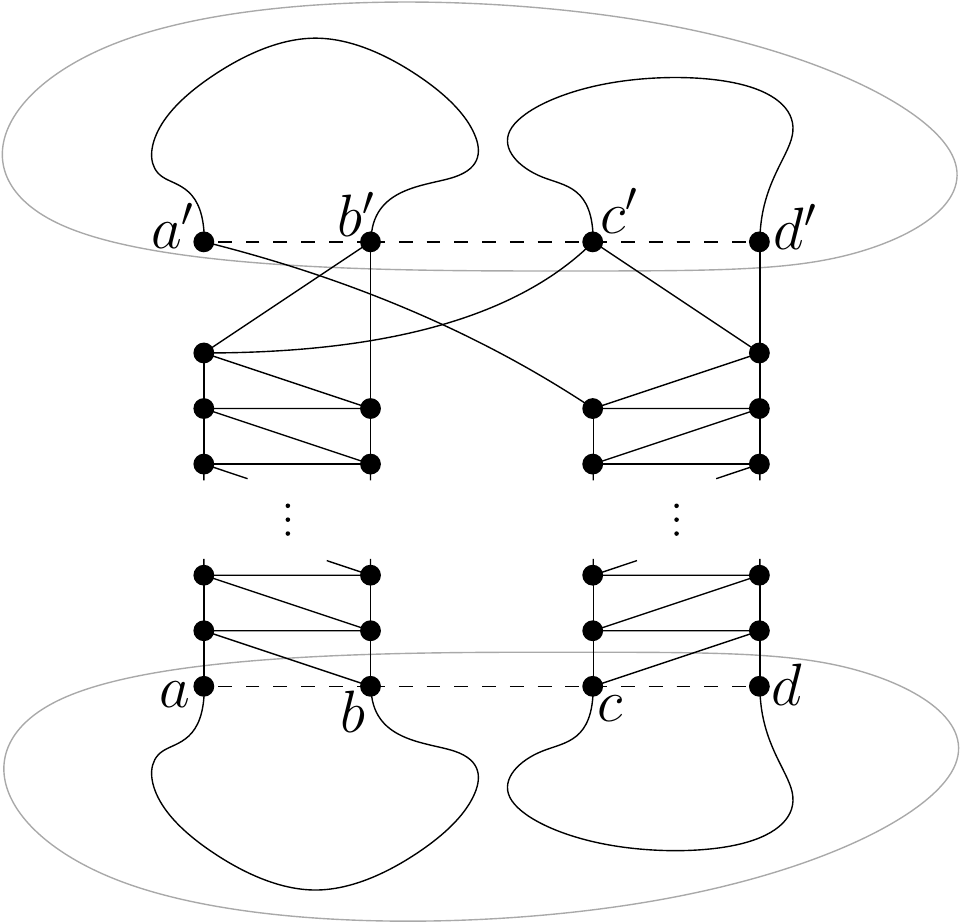}\includegraphics[height=65mm]{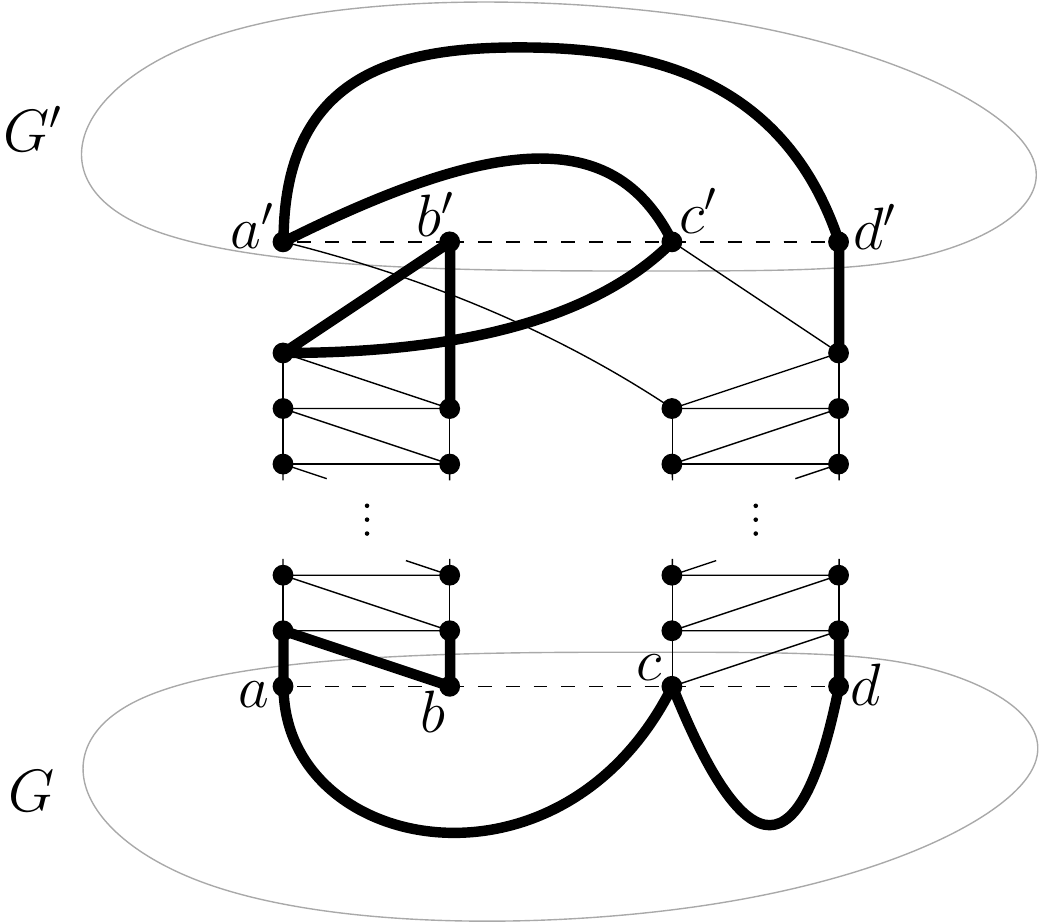}\\
Fig.~2: On the left, the graph $H_G$ is shown, with relevant 2-factors in $G$ and $G'$. On the right, we illustrate why conditions (i), (ii), and the first statement in (iii) would not suffice: if $G$ and $G'$ are traversed as shown, then the number of Hamiltonian cycles in $H_G$ would increase with $\ell$.
\end{center}

The needed graph $G$ shall be the modification of the Petersen graph shown in Figure 3 where the bold edges form the 2-factor satisfying (ii). It is straightforward to infer from the non-Hamiltonicity of Petersen's graph that $G$ satisfies (i) and (iii).

\begin{center}
\includegraphics[height=40mm]{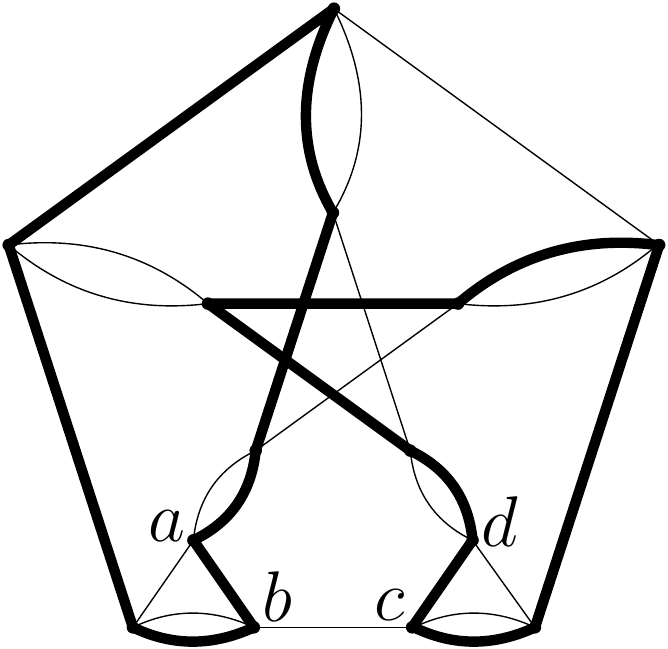}\\
Fig.~3: A useful modification of Petersen's graph.
\end{center}

Using two copies of $G$, we construct $H_G$ as described above.
One final issue remains: there are double edges occurring in $H_G$. We shall make use of the idea behind Meredith's classical construction in which a vertex is replaced by a complete bipartite graph $K_{3,4}$, see~\cite[p.~161]{BM76}. In the Meredith graph the operation is performed on each vertex of the Petersen graph (in which a 1-factor is replaced by double edges). In the present note the operation is performed only on both ends of each double edge in $H_G$. \hfill $\Box$

\vspace{3cm}

\noindent \textbf{Acknowledgements.} Thomassen's research was supported by the Independent Research Fund Denmark, 8021-002498 AlgoGraph. Zamfirescu's research was supported by a Postdoctoral Fellowship of the Research Foundation Flanders (FWO).

\end{document}